\documentclass[12pt]{amsart}
\usepackage{amsmath,amssymb,graphicx}
\usepackage{hyperref}

\newtheorem{theorem}{Theorem}
\newtheorem{lemma}[theorem]{Lemma}

\newtheorem{conjecture}{Conjecture}
\newtheorem*{question}{Question}

\theoremstyle{remark}
\newtheorem{remark}{Remark}

\newcommand{\bbK}{{\mathbb K}}
\newcommand{\bbQ}{{\mathbb Q}}
\newcommand{\bbR}{{\mathbb R}}
\newcommand{\bbZ}{{\mathbb Z}}

\newcommand{\ba}{{\mathbf a}}
\newcommand{\balpha}{\boldsymbol{\alpha}}
\newcommand{\lcm}{{\rm lcm}}

\begin{document}

\title[Families of periodic Jacobi-Perron algorithms]{Families of periodic Jacobi-Perron algorithms for all period lengths}

\author{Paul Voutier}
\curraddr{London, UK}
\email{paul.voutier@gmail.com}

\begin{abstract}
For all integers $m \geq  n \geq 2$, we exhibit infinite families of purely
periodic Jacobi-Perron Algorithm (JPA) expansions of dimension $n$ with period length
equal to $m$ along with the associated Hasse-Bernstein units.
Some observations on the units of Levesque-Rhin as well as the periodicity
of the JPA expansion of $\left( m^{1/3}, m^{2/3} \right)$ are also made.
\end{abstract}

\maketitle

\section{Introduction}

The Jacobi-Perron Algorithm (JPA) is one of the generalisations of the
ordinary continued fraction algorithm to higher dimensions. It has been the
subject of considerable study (see, for example, \cite{B3,S}).

Let $\balpha^{(0)}=\left( \alpha_{1}^{(0)},\ldots, \alpha_{n-1}^{(0)} \right)
\in \bbR^{n-1}$ for $n \geq 2$. The JPA expansion of $\balpha^{(0)}$
is the sequence of elements of $\bbR^{n-1}$, $\langle \balpha^{(\nu)} \rangle_{\nu \geq 0}$,
defined by
$$
\begin{array}{ll}
\ba^{(\nu)}       &= \left( a_{1}^{(\nu)},                      \ldots, a_{n-1}^{(\nu)} \right) \\[3.0mm]
                  &= \left( \lfloor \alpha_{1}^{(\nu)} \rfloor, \ldots, \lfloor \alpha_{n-1}^{(\nu)} \rfloor \right), \\[3.0mm]
\balpha^{(\nu+1)} &= \left( \alpha_{1}^{(\nu+1)},               \ldots, \alpha_{n-2}^{(\nu+1)}, \alpha_{n-1}^{(\nu+1)} \right) \\[3.0mm]
                  &= \left( \displaystyle\frac{\alpha_{2}^{(\nu)}-a_{2}^{(\nu)}}{\alpha_{1}^{(\nu)}-a_{1}^{(\nu)}},
                            \ldots, \displaystyle\frac{\alpha_{n-1}^{(\nu)}-a_{n-1}^{(\nu)}}{\alpha_{1}^{(\nu)}-a_{1}^{(\nu)}},
                            \displaystyle\frac{1}{\alpha_{1}^{(\nu)}-a_{1}^{(\nu)}} \right),
\end{array}
$$
where $\alpha_{1}^{(\nu)} \neq a_{1}^{(\nu)}$ and $\lfloor \cdot \rfloor$ is the greatest
integer function.

For $n=2$, this is the ordinary continued fraction algorithm.
When $n=3$, this algorithm was given by Jacobi \cite{J} and is sometimes
called the {\it Jacobi Algorithm}. Perron \cite{P} extended Jacobi's
definition to arbitrary $n$. As a result, for all $n \geq 2$, the algorithm
is now known as the {\it Jacobi-Perron Algorithm}. Following \cite{B3,LR}, we
have presented it here in its inhomogeneous form.

The JPA expansion of $\balpha^{(0)}$ is called {\it periodic}, if there exist
two integers $\ell_{0}$, $\ell_{1}$ with $\ell_{0} \geq 0$ and $\ell_{1} \geq 1$
such that
$$
\balpha^{(\nu+\ell_{1})}=\balpha^{(\nu)}
\text{ for } \nu=\ell_{0},\ell_{0}+1,\ldots .
$$

If $\ell_{0}$ and $\ell_{1}$ are the smallest integers satisfying these
conditions, then
$$
\balpha^{(0)}, \balpha^{(1)},\ldots, \balpha^{(\ell_{0}-1)}
\text{ and }
\balpha^{(\ell_{0})}, \balpha^{(\ell_{0}+1)},\ldots, \balpha^{(\ell_{0}+\ell_{1}-1)}
$$
are called respectively the {\it preperiod} and the {\it period} of the periodic
JPA expansion, and $\ell_{0}$ and $\ell_{1}$ are their respective {\it lengths}.
If $\ell_{0}=0$, then the JPA expansion of $\balpha^{(0)}$ is said to be
{\it purely periodic}.

When the JPA expansion is periodic and all the $\alpha_{i}^{(0)}$'s are algebraic
integers, then Hasse and Bernstein \cite{BH} proved that
$$
\epsilon = \prod_{i=\ell_{0}}^{\ell_{0}+\ell_{1}-1} \alpha_{n-1}^{(i)}
$$
is a unit in the ring of algebraic integers of
$\bbQ \left( \alpha_{1}^{(0)},\alpha_{2}^{(0)},\ldots,\alpha_{n-1}^{(0)} \right)$.
This is the {\it Hasse-Bernstein formula} and $\epsilon$ is called the
{\it Hasse-Bernstein unit}.

Lagrange showed that the ordinary continued fraction expansion of a real number
$\alpha$ is periodic if and only if $\alpha$ is a quadratic irrational.
A central problem in the study of multidimensional generalisations
of the ordinary continued fraction expansion is to find, and prove, an
analogue of this result. The following question is part of JPA folklore.

\begin{question}
When $1,\alpha_{1},\ldots,\alpha_{n-1} \in \bbR$ form a $\bbQ$-basis of a number field
of degree $n$, is the JPA expansion of $\left( \alpha_{1}, \ldots,
\alpha_{n-1} \right)$ periodic?
\end{question}

Numerical evidence suggests this is not always true for $n \geq 3$. There
are many such tuples for which we do not know if the JPA expansion is periodic.
For example, with $\left( \alpha_{1}, \alpha_{2} \right)
= \left( m^{1/3}, m^{2/3} \right)$ our own calculations (see \cite{V1}) show
that $\ell_{0}+\ell_{1}>30,000$ for many $4 \leq m \leq 5000$, including $m=4$,
$6,11,13$, $15,19,20$. In contrast, for all such JPA expansions known to be
periodic with $m \leq 5000$, the maximum value of $\ell_{0}+\ell_{1}$ is $93$,
occurring for $m=17$.

Furthermore, also in the forthcoming \cite{V1}, we present results suggesting
structure to the periodic behaviour of such pure cubic (and $n$-th root) JPA
expansions as well as supporting evidence for the following conjecture.

\begin{conjecture}
Let $m=x^{3}-x$ where $x \geq 2$ is a positive integer. The JPA expansion of
$\left( m^{1/3}, m^{2/3} \right)$ is not periodic.
\end{conjecture}

$x^{3}-x$ appears to be one of many expressions of the form $ax^{3}+bx$ with
$a,b \in \bbZ$ for which this conjecture is true.

When $n=3$ (i.e., the cubic case), only a few infinite classes of periodic
Jacobi Algorithm expansions were known and only for period length at most 9,
until 1991 when Levesque and Rhin \cite{LR} produced infinite classes for all
period lengths of the form $3m+1$ and $4m+1$. They also determined the associated
Hasse-Bernstein units, $\omega$. Adam \cite{Adam2} proved that these $\omega$
are fundamental units in the ring of integers, $\bbZ [ \omega ]$.

In Section~\ref{subsect:units-LR}, we provide examples answering negatively
the questions of Levesque and Rhin \cite{LR} on whether their Hasse-Bernstein
units are fundamental units in the full ring of integers of $\bbQ(\omega)$.

For $n \geq 4$, less is known. Bernstein \cite{B1,B2} (see also \cite{L,DP3})
showed that one can construct purely periodic JPA expansions of period lengths
$1$, $n$ and $n^{2}$ from $\left( D^{n} \pm d \right)^{1/n}$ when $n \geq 3$,
$d|D$ and $D$ is sufficiently large relative to $d$. In 1984, Paysant-Le Roux
and Dubois \cite{DP4} showed that every real number field, $\bbK$, of degree
$n+1$ contains infinitely many $n$-tuples $\left( \alpha_{1}, \ldots, \alpha_{n} \right)$
such that the JPA expansion of $\left( \alpha_{1}, \ldots, \alpha_{n} \right)$
is purely periodic of length $n+1$ and $\left\{ 1, \alpha_{1},\ldots,\alpha_{n} \right\}$
is a $\bbQ$-basis for $\bbK$. Also, in her thesis, Adam \cite{Adam1} generalised
the results of Levesque and Rhin to all $n \geq 2$ and all periods of the
form $kn+1$, where $n \leq k \leq 2n-2$.

\section{Results}

Here we extend the results of Levesque and Rhin to all dimensions, $n \geq 2$,
and all period lengths, $m \geq n$, producing infinite families of expansions
for each such $m$ and $n$. In our expansions, the terms grow linearly for a
fixed period length, $m$, rather than polynomially of degree $m$ as in their work.

We start by defining $n$ related recurrence sequences, $\left\{ u_{0,m} \right\}_{m \geq 0}$,
$\ldots$,$\left\{ u_{n-1,m} \right\}_{m \geq 0}$, for any $n$-tuple
$\left( c_{0}, c_{1}, \ldots, c_{n-1} \right)$ of non-negative integers with
$c_{0},c_{n-1} \geq 1$. Set $u_{i,m}=\delta_{i,m}$, the Kronecker delta function,
for $0 \leq i, m \leq n-1$, and define
\begin{equation}
\label{eq:recur}
u_{i,m}=c_{n-1}u_{i,m-1}+c_{n-2}u_{i,m-2}+\cdots+c_{0}u_{i,m-n}
\end{equation}
for $i=0,1,\ldots,n-1$ and all $m \geq n$.

\begin{theorem}
\label{thm:1}
For all integers $c_{0},c_{1}, \ldots, c_{n-1}, m,n,t$ with $c_{n-1} \geq c_{0}=1$,
$c_{n-1} \geq c_{n-2} \geq \cdots \geq c_{1} \geq 0$, $m \geq n \geq 2$, and
$t \geq 1$, let
\[
f(X)= X^{n} - \left( u_{n-1,m}t+c_{n-1} \right) X^{n-1} - \cdots - \left( u_{0,m}t+c_{0} \right)
\]
and let $\omega$ be the unique positive real root of $f(X)$.
The JPA expansion of $\balpha^{(0)}$ with
\begin{equation}
\label{eq:exp1}
\alpha_{i}^{(0)}
=\sum_{j=0}^{i} \left( u_{j,m}t+c_{j} \right) \omega^{j-i}
\end{equation}
for $i=1,\ldots,n-1$, is purely periodic of period length $m$
and its $\ba^{(\nu)}$'s are as in Table~$\ref{table:a-values}$.
\begin{table}[h]
\scalebox{0.9}{%
\begin{tabular}{ll}
\hline
$\nu \bmod m$       & $\ba^{(\nu)}$                                               \\ \hline
$0$                 & $\left( u_{1,m}t+c_{1}, \ldots, u_{n-1,m}t+c_{n-1} \right)$ \\
$1,\ldots,m-n+1$    & $\left(          c_{1}, \ldots,            c_{n-1} \right)$ \\
$m-n+2,\ldots,m-1$  & $\left(          c_{1}, \ldots, c_{m-\nu}, u_{m-\nu+1,m}t+c_{m-\nu+1}, \ldots, u_{n-1,m}t+c_{n-1} \right)$ \\ \hline
\end{tabular}%
}
\label{table:a-values}
\caption{$\ba^{(\nu)}$'s}
\end{table}

Furthermore, the Hasse-Bernstein unit associated with this JPA expansion is
\[
\epsilon= u_{n-1,m}\omega^{n-1}+u_{n-2,m}\omega^{n-2}+\cdots+u_{1,m}\omega+u_{0,m}.
\]
\end{theorem}

\begin{remark}
\label{rem:1}
Applying $f(\omega)=0$ to the expression in \eqref{eq:exp1}, we can also write
$\balpha^{(0)}$ as
$$
\left( \omega^{n-1}- \sum_{i=1}^{n-2} \left( u_{i+1,m}t+c_{i+1} \right) \omega^{i},
\ldots, \omega^{2}-\left( u_{n-1,m}t+c_{n-1} \right) \omega, \omega \right).
$$
This choice of our starting vector follows the advice and experience of Dubois
and Paysant-Le Roux with finding periodic JPA expansions \cite{DP1,DP2,DP3},
and in particular purely periodic ones.
\end{remark}

\begin{remark}
If $c_{0}>1$, then we get a purely periodic expansion of period length $\lcm (m,n)$
provided that $c_{0}|c_{i}$ for all $i=1,\ldots,n-1$. It consists of $\lcm(m,n)/m$
repetitions of the above pattern of $\ba^{(\nu)}$'s but with $c_{0}$
in some of the denominators.

If $c_{0}$ does not divide all the other $c_{i}$'s, then computational evidence
suggests that the expansion is not periodic.
\end{remark}

\section{Preliminary Lemmas}

\begin{lemma}
\label{lem:3}
Let $n \geq 2$ be an integer and
$$
g(X)=X^{n}-a_{n-1}X^{n-1}-a_{n-2}X^{n-2}-\cdots-a_{0} \in \bbR[X]
$$
with $a_{n-1} \geq 1$, $a_{n-1} \geq a_{0}>0$ and
$a_{n-1} \geq a_{n-2} \geq \cdots \geq a_{1} \geq 0$. Then $g(X)$ has exactly
one non-negative real root, $\omega$; it satisfies $a_{n-1}<\omega<a_{n-1}+1$.
For all $1 \leq i \leq n-1$, we have
\begin{equation}
\label{eq:lem3}
0<\sum_{j=0}^{i-1} a_{j}\omega^{j-i}<1.
\end{equation}
\end{lemma}

\begin{proof}
We first establish the statement about the non-negative real root of $g(X)$
and its location.

If $0 \leq x \leq a_{n-1}$, then $g(x) \leq -a_{0}<0$.
From our condition on the coefficients of $g(X)$,
\begin{eqnarray*}
\frac{g \left( a_{n-1}+1 \right)}{\left( a_{n-1}+1 \right)^{n-2}}
& \geq & \left( a_{n-1}+1 \right)^{2}-a_{n-1}\left( a_{n-1}+1 \right)-
         a_{n-1} \sum_{i=0}^{n-2} \left( a_{n-1}+1 \right)^{-i} \\
& > & \left( a_{n-1}+1 \right)^{2}-a_{n-1}\left( a_{n-1}+1 \right)-2a_{n-1}
=a_{n-1}-1 \geq 0,
\end{eqnarray*}
since $a_{n-1} \geq 1$. Therefore, $g \left( a_{n-1}+1 \right) \geq 0$.

As $g'(x)>n \left( g(x)+a_{0} \right)/x$ for $x>0$, so $g'(x)>0$ for
$x \geq \omega$. Hence there is only one non-negative real root of $g(X)$.

Since $a_{0}, \omega>0$, the lower bound in \eqref{eq:lem3} holds.

The upper bound in \eqref{eq:lem3} is equivalent to
$\sum_{j=0}^{i-1}a_{j}\omega^{j}<\omega^{i}$.
We know from $g(\omega)=0$, $a_{0}>0$, $a_{n-1} \geq \cdots \geq a_{1}$ and
$i \leq n-1$, that
$$
\omega^{n}>\sum_{j=0}^{i-1}a_{n+j-i}\omega^{n+j-i}
\geq \sum_{j=0}^{i-1}a_{j}\omega^{n+j-i}.
$$
Dividing by $\omega^{n}$ completes the proof.
\end{proof}

We also need relations between the elements of the recurrence sequences defined
at the start of Section~2. Let $c_{0}, c_{1}, \ldots, c_{n-1}$ and
$\left\{ u_{0,m} \right\}_{m \geq 0}$, $\ldots$,
$\left\{ u_{n-1,m} \right\}_{m \geq 0}$ be as there.

\begin{lemma}
\label{lem:2a}
{\rm (a)} For $m \geq 1$,
\begin{equation}
\label{eq:1}
u_{0,m}=c_{0}u_{n-1,m-1}
\end{equation}
and for $i=1,\ldots,n-1$,
\begin{equation}
\label{eq:2}
u_{i,m}=c_{i}u_{n-1,m-1}+u_{i-1,m-1}.
\end{equation}

\noindent
{\rm (b)} For $i=1,\ldots,n-1$ and $m \geq i+1$,
$$
u_{i,m}=\sum_{j=0}^{i} c_{j}u_{n-1,m-i+j-1}.
$$
\end{lemma}

\begin{proof}
(a) These relationships are immediately verified for $m \leq n$ using the initial
conditions in the definitions of the sequences and observing that $u_{i,n}=c_{i}$.
Applying the recurrence relation \eqref{eq:recur} and induction completes the
proof for $m>n$.

(b) We prove this relation using part~(a) and induction.

For $i=1$,
$$
u_{1,m}=c_{1}u_{n-1,m-1}+u_{0,m-1}=c_{1}u_{n-1,m-1}+c_{0}u_{n-1,m-2},
$$
as required.

Now consider $2 \leq i \leq n-1$. From \eqref{eq:2}, we have
$u_{i,m}=c_{i}u_{n-1,m-1}+u_{i-1,m-1}$ and the desired relation now follows
immediately by induction.
\end{proof}

Finally, we need some inequalities satisfied by elements of these recurrence
sequences too.

\begin{lemma}
\label{lem:2b}
\noindent
{\rm (a)} For $m \geq n$, $u_{i,m} \geq u_{i,m-1}$ for $i=0,\ldots,n-1$ and
for at least one such $i$, $u_{i,m}>u_{i,m-1}$.

\noindent
{\rm (b)} For $m \geq 0$, $u_{i,m} \geq 1$ for at least one $i=0,\ldots,n-1$.

\noindent
{\rm (c)} Let $m \geq n-1$. Then $u_{i,m} \geq u_{i-1,m}$ for $2 \leq i \leq n-1$
and $u_{n-1,m} \geq u_{0,m}$.

\noindent
{\rm (d)} Suppose that $c_{0}=1$ and $c_{n-1} \geq c_{n-2} \geq \cdots \geq c_{1}$.
For $m \geq n$ and $1 \leq j \leq n-1$,
\begin{equation}
\label{eq:2b-d}
c_{j-1}u_{i,m-1}+\cdots+c_{1}u_{i,m-j+1}+c_{0}u_{i,m-j} \leq u_{i,m}+1
\end{equation}
for all $0 \leq i \leq n-1$. For each $1 \leq j \leq n-1$, there exists
$i_{1} \in \left\{ 1, \ldots, n-1 \right\}$ such that
$c_{j-1}u_{i_{1},m-1}+\cdots+c_{1}u_{i_{1},m-j+1}+c_{0}u_{i_{1},m-j}<u_{i_{1},m}$,
and there exists such an $i_{1}$ larger than any $i$ for which equality holds
in \eqref{eq:2b-d}.
\end{lemma}

\begin{proof}
(a) The first statement follows by the recurrence relation \eqref{eq:recur},
since $c_{n-1} \geq 1$.

Since $c_{0} \geq 1$, for $n \leq m \leq 2n-1$, we have
$$
u_{m-n,m} \geq c_{n-1}u_{m-n,m-1}+c_{0}u_{m-n,m-n}
\geq u_{m-n,m-1}+1.
$$

For $m \geq 2n$, we can take $i=0$ and apply the recurrence relation \eqref{eq:recur},
since $u_{0,n},\ldots,u_{0,2n-1} \geq 1$. This establishes the second statement
in part~(a).

(b) For $m \geq n$, this is immediate from (a). For $0 \leq m \leq n-1$, it
follows from the initial values for the $u_{i,m}$'s.

(c) Applying first \eqref{eq:2} with $i=n-1$ and then \eqref{eq:1} from
Lemma~\ref{lem:2a}(a), we have $u_{n-1,m}=\left( c_{n-1}/c_{0} \right)u_{0,m}+u_{n-2,m-1}
\geq u_{0,m}$ by our assumption that $c_{n-1} \geq c_{0}$ and since the $u_{i,j}$'s
are non-negative. So we need only consider the first part of (c) in what follows.

This is true for $m=n-1$, since $u_{n-1,n-1}=1$ and $u_{i,n-1}=0$ for
$i=0,\ldots,n-2$.

Let $m$ be an integer with $n \leq m \leq 2n-1$ and suppose that (c) holds
for $n-1 \leq j < m$. We can write
$$
u_{i,m}=
\left\{
\begin{array}{ll}
c_{n-1}u_{i,m-1}+\cdots+c_{2n-m}u_{i,n}+c_{n-m+i} & \mbox{if $i \geq m-n$} \\
c_{n-1}u_{i,m-1}+\cdots+c_{2n-m}u_{i,n} & \mbox{otherwise},
\end{array}
\right.
$$
since $c_{2n-m}u_{i,n}=c_{n-(m-n)}u_{i,m-(m-n)}$ and $u_{j,k}=\delta_{j,k}$
for $0 \leq j,k \leq n-1$, by the definition of these recurrence sequences.

Part~(c) now follows for $m$ by our assumption that it holds for
$n-1 \leq j < m$ and our assumption on the sizes of the $c_{i}$'s.

For $m \geq 2n$, our desired result follows immediately from the fact that it
holds for $n-1 \leq m \leq 2n-1$ and the recurrence relation \eqref{eq:recur}
defining the sequences.

(d) We start by showing that the stronger inequality
\begin{equation}
\label{eq:strong-e}
c_{j-1}u_{i,m-1}+\cdots+c_{1}u_{i,m-j+1}+u_{i,m-j}
\leq u_{i,m},
\end{equation}
holds for $m \geq n$ if $c_{n-j} \geq 1$ and for $m-j \geq n-1$ otherwise (this
condition is best possible: with $c_{0}=c_{n-1}=1$ and $c_{1}=\cdots=c_{n-2}=0$,
$i=n-2$, $j=n-1$ and $m=2n-3$, we find that $u_{i,m}=0$, while
$c_{j-1}u_{i,m-1}+\cdots+c_{1}u_{i,m-j+1}+u_{i,m-j}=1$).

For $m \geq n$, using the recurrence relation \eqref{eq:recur}, to establish
\eqref{eq:strong-e} we need to show that
\begin{eqnarray}
\label{eq:lem4-d}
& & c_{j-1}u_{i,m-1}+\cdots+c_{1}u_{i,m-j+1}+u_{i,m-j} \\
& \leq &
c_{n-1}u_{i,m-1}+\cdots+c_{n-j+1}u_{i,m-j+1}+c_{n-j}u_{i,m-j}+\cdots+u_{i,m-n}, \nonumber
\end{eqnarray}
for all $1 \leq i \leq n-1$. Since $c_{1} \leq \cdots \leq c_{n-1}$, this will
follow provided $c_{n-j} \geq 1$ (or if $u_{i,m-j}=0$, but we will not use that
here).

Suppose that $c_{n-j}=0$ and let $k$ be the largest positive integer less than
$j$ such that $c_{n-k} \geq 1$. Since $j \leq n-1$, we have $c_{j-k}=0$ and
provided that $u_{i,m-k} \geq u_{i,m-j}$, \eqref{eq:strong-e} follows. Moreover,
$u_{i,m-k} \geq u_{i,m-j}$ holds for $m-j \geq n-1$ by part~(a).

Thus \eqref{eq:strong-e} holds under the stated conditions that follow it.

We now turn to proving that \eqref{eq:2b-d} holds. We need only consider the
case when $m-j<n-1$ and $c_{n-j}=0$. But in this case, $u_{i,m-j} \leq 1$, while
the right-hand side of \eqref{eq:lem4-d} is greater than or equal to the sum
of the first $j-1$ terms on the left-hand side. Hence \eqref{eq:2b-d} follows.

Next, we show the existence of $i_{1}$.

By part~(b), $u_{i,m-n} \geq 1$ for at least one $i$ for $m \geq n$, so under
the conditions above (i.e., $m \geq n$ if $c_{n-j} \geq 1$ and $m-j \geq n-1$
otherwise), there is at least one $i$ for which we have strict inequality in
\eqref{eq:strong-e} too. The case when $m-j < n-1$ will be covered in the next
paragraph.

It only remains to show that if we have equality in \eqref{eq:2b-d} for $i$,
then we have strict inequality in \eqref{eq:strong-e} for some $i_{1}>i$. From
our analysis above for when \eqref{eq:strong-e} holds, we know that we can only
have equality in \eqref{eq:2b-d} if $m-j<n-1$ and $c_{n-j}=0$. So this can only
occur if $u_{i,m-j}=1$ (i.e., $i=m-j$). We take
$i_{1}=\min \left( m-k, n-1 \right)>m-j=i$, so
\begin{eqnarray*}
& & c_{j-1}u_{i_{1},m-1}+\cdots+c_{j-k}u_{i_{1},m-k}+\cdots+c_{1}u_{i_{1},m-j+1}+u_{i_{1},m-j} \\
& \leq &
c_{n-1}u_{i_{1},m-1}+\cdots+c_{n-k}u_{i_{1},m-k}=u_{i_{1},m},
\end{eqnarray*}
since $u_{i_{1},m-j}=0$, $c_{n-1} \geq \cdots \geq c_{1}$ and $c_{n-k-1}=\cdots=c_{1}=0$.

Notice that
$u_{i_{1},m-k} \geq 1$, by our choice of $i_{1}$. By our choice of $k$, we have
$c_{n-k}=1$ and $c_{j-k}=0$. Thus $c_{n-k}u_{i_{1},m-k} \geq 1$ and $c_{j-k}u_{i_{1},m-k}=0$.
Therefore, we have strict inequality and for a value of $i_{1}$ satisfying the
required conditions.
\end{proof}

\section{Proof of Theorem~\ref{thm:1}}

\subsection{$\ba^{(0)}$}

From \eqref{eq:exp1}, we have
$\alpha_{i}^{(0)}- \left( u_{i,m}t+c_{i} \right)
=\sum_{j=0}^{i-1} \left( u_{j,m}t+c_{j} \right) \omega^{j-i}$.
From Lemma~\ref{lem:2b}(c) and our assumptions in Theorem~\ref{thm:1} on the
size of the $c_{j}$'s, we can apply Lemma~\ref{lem:3} with $g(X)=f(X)$ and
$a_{i}=u_{i,m}t+c_{i}$. From \eqref{eq:1} in Lemma~\ref{lem:3}, the above sum for
$\alpha_{i}^{(0)}- \left( u_{i,m}t+c_{i} \right)$ is strictly between $0$ and
$1$ and so
\begin{equation}
\label{eq:exp-ba}
\ba^{(0)} = \left( u_{1,m}t+c_{1}, \ldots, u_{n-1,m}t+c_{n-1} \right).
\end{equation}

\subsection{Expression for $\balpha^{(1)}$ independent of $t$}

A crucial fact for us about $\balpha^{(1)}$ is that it can be written independently
of $t$. We establish such an expression now.

From \eqref{eq:exp1} and the above expression in \eqref{eq:exp-ba} for $\ba^{(0)}$,
we obtain
\begin{equation}
\label{eq:exp2}
\alpha_{i}^{(0)}-a_{i}^{(0)}
= \omega^{-i} \sum_{j=0}^{i-1} \left( u_{j,m}t+c_{j} \right) \omega^{j}.
\end{equation}

For $i=1,\ldots,n-1$, we consider
$$
\left( \alpha_{i}^{(0)}-a_{i}^{(0)} \right)
\left( u_{n-1,m}\omega^{n-1}+\cdots+u_{0,m} \right).
$$

Applying \eqref{eq:exp2}, we have
\begin{eqnarray*}
&   & \left( \alpha_{i}^{(0)}-a_{i}^{(0)} \right)
      \left( u_{n-1,m}\omega^{n-1}+\cdots+u_{0,m} \right) \\
& = & \omega^{-i} \left( \sum_{j=0}^{i-1} \left( u_{j,m}t+c_{j} \right) \omega^{j} \right)
      \left( u_{n-1,m}\omega^{n-1}+\cdots+u_{0,m} \right) \\
& = & \omega^{-i} \left\{
      t\left( u_{n-1,m}\omega^{n-1}+\cdots+u_{0,m} \right)
      \left( \sum_{j=0}^{i-1} u_{j,m} \omega^{j} \right) \right. \\
&   & \left. +\left( u_{n-1,m}\omega^{n-1}+\cdots+u_{0,m} \right)
\left( \sum_{j=0}^{i-1} c_{j} \omega^{j} \right) \right\}.
\end{eqnarray*}

Since $f(\omega)=0$, we can write this last expression as
\begin{eqnarray*}
&   & \omega^{-i} \left\{
      \left( \omega^{n}-c_{n-1}\omega^{n-1}-\cdots-c_{i}\omega^{i}
      -\left( \sum_{j=0}^{i-1} c_{j}\omega^{j} \right) \right)
      \left( \sum_{j=0}^{i-1} u_{j,m} \omega^{j} \right) \right. \\
&   & \left. +\left( u_{n-1,m}\omega^{n-1}+\cdots+u_{i,m}\omega^{i}
      + \sum_{j=0}^{i-1} u_{j,m} \omega^{j} \right)
      \left( \sum_{j=0}^{i-1} c_{j} \omega^{j} \right) \right\} \\
& = & \left( \omega^{n-i}-\sum_{j=i}^{n-1}c_{j}\omega^{j-i} \right)
      \left( \sum_{j=0}^{i-1} u_{j,m} \omega^{j} \right)
      +\left( \sum_{j=i}^{n-1} u_{j,m}\omega^{j-i} \right)
      \left( \sum_{j=0}^{i-1} c_{j} \omega^{j} \right).
\end{eqnarray*}

We now show that for $i=1,\ldots,n-1$,
\begin{eqnarray}
\label{eq:alphai-simp}
& & \left( \omega^{n-i}-\sum_{j=i}^{n-1}c_{j}\omega^{j-i} \right)
      \left( \sum_{j=0}^{i-1} u_{j,m} \omega^{j} \right)
      +\left( \sum_{j=i}^{n-1} u_{j,m}\omega^{j-i} \right)
      \left( \sum_{j=0}^{i-1} c_{j} \omega^{j} \right) \\
& = & \sum_{k=0}^{n-1} \left( \sum_{j=0}^{i-1} c_{j}u_{k,m-i+j} \right) \omega^{k}.
\nonumber
\end{eqnarray}

From this, along with the relationship between $\balpha^{(0)}$, $\ba^{(0)}$ and
$\balpha^{(1)}$ that defines a JPA, it follows that for $i=1,\ldots,n-1$,
\begin{equation}
\label{eq:ai}
\alpha_{i}^{(1)} = \frac{\sum_{k=0}^{n-1} \left( \sum_{j=0}^{i} c_{j}u_{k,m-1-i+j} \right) \omega^{k}}
{c_{0} \left( u_{n-1,m-1}\omega^{n-1}+\cdots+u_{0,m-1} \right)}.
\end{equation}
This is immediate for $i=1,\ldots,n-2$. For $i=n-1$, we just need to apply
\eqref{eq:recur} to the coefficient of each $\omega^{k}$ term in the numerator
to establish the equality.

We prove our simplification by induction. 
Applying Lemma~\ref{lem:2a}(a), for $i=1$, the left-hand side of
\eqref{eq:alphai-simp} is $u_{n-1,m-1}\omega^{n-1}+\cdots+u_{0,m-1}$, as
required.

We now assume that \eqref{eq:alphai-simp} holds for some $1 \leq i \leq n-1$
and use some notation to assist our argument. Define
\begin{eqnarray*}
I_{1}^{(i)} & = & \omega^{n-i}-\sum_{j=i}^{n-1}c_{j}\omega^{j-i},
\hspace{13.0mm}
I_{2}^{(i)} =  \sum_{j=0}^{i-1} u_{j,m} \omega^{j}, \\
I_{3}^{(i)} & = & \sum_{j=i}^{n-1} u_{j,m}\omega^{j-i}
\hspace{10.0mm} \mbox{ and } \hspace{3.0mm}
I_{4}^{(i)} = \sum_{j=0}^{i-1} c_{j} \omega^{j}.
\end{eqnarray*}

The left-hand side of \eqref{eq:alphai-simp} can be written as
$I_{1}^{(i)}I_{2}^{(i)}+I_{3}^{(i)}I_{4}^{(i)}$ and we have
\begin{eqnarray*}
I_{1}^{(i+1)} & = & \frac{I_{1}^{(i)}+c_{i}}{\omega}, \hspace{15.5mm}
I_{2}^{(i+1)}=I_{2}^{(i)}+u_{i,m}\omega^{i}, \\
I_{3}^{(i+1)} & = & \frac{I_{3}^{(i)}-u_{i,m}}{\omega}
\hspace{2.0mm} \mbox{ and } \hspace{2.0mm}
I_{4}^{(i+1)}=I_{4}^{(i)}+c_{i}\omega^{i}.
\end{eqnarray*}

Expanding these expressions, we find that
\begin{eqnarray*}
&   & \omega \left( I_{1}^{(i+1)}I_{2}^{(i+1)}+I_{3}^{(i+1)}I_{4}^{(i+1)} \right) \\
& = & \left( I_{1}^{(i)}I_{2}^{(i)}+I_{3}^{(i)}I_{4}^{(i)} \right)
+u_{i,m}\omega^{n}+\sum_{k=0}^{n-1} \left( c_{i}u_{k,m}-c_{k}u_{i,m} \right) \omega^{k}.
\end{eqnarray*}

The desired expression for the left-hand side of this equation can be written as
\begin{equation}
\label{eq:lhs}
\left( \sum_{j=0}^{i} c_{j}u_{n-1,m-i+j-1} \right) \omega^{n}
+ \sum_{k=1}^{n-1} \left( \sum_{j=0}^{i} c_{j}u_{k-1,m-i+j-1} \right) \omega^{k},
\end{equation}
while the right-hand side is
\begin{equation}
\label{eq:rhs}
u_{i,m}\omega^{n}
+ \sum_{k=1}^{n-1} \left( -c_{k}u_{i,m}+\sum_{j=0}^{i} c_{j}u_{k,m-i+j-1} \right) \omega^{k}
-c_{0}u_{i,m}+\sum_{j=0}^{i}c_{j}u_{0,m-i+j}.
\end{equation}

Using Lemma~\ref{lem:2a}(b), we see that the coefficients of $\omega^{n}$ are
equal. Similarly, using \eqref{eq:1} from Lemma~\ref{lem:2a}(a) along with
Lemma~\ref{lem:2a}(b), the constant coefficient in \eqref{eq:rhs} is zero.

Finally, we compare the coefficients for $k=1,\ldots,n-1$. Using \eqref{eq:1},
we find that the coefficient of $\omega^{k}$ in \eqref{eq:rhs} minus its
counterpart in \eqref{eq:lhs} is
$$
-c_{k}u_{i,m}+c_{k}\sum_{j=0}^{i} c_{j}u_{n-1,m-i+j-1},
$$
which is zero by again applying Lemma~\ref{lem:2a}.

This completes our proof of \eqref{eq:alphai-simp} and thus \eqref{eq:ai}.

\subsection{Induction from $\balpha^{(\nu)}$ to $\balpha^{(\nu+1)}$:
$1 \leq \nu \leq m-n$}
$\,$

At this point, we use our assumption that $c_{0}=1$, the reason being that
some of the $\alpha_{i}^{(\nu)}$'s below are divided by $c_{0}$, depending on
the values of both $i$ and $\nu$. This complicates the expressions in this
subsection and the next one as well as resulting in a longer period, as stated
above after Theorem~\ref{thm:1}.

Suppose that for some $\nu$ satisfying $1 \leq \nu \leq m-n$ and for all $i=1,\ldots,n-1$,
we have
$$
\alpha_{i}^{(\nu)} = \frac{\sum_{k=0}^{n-1} \left( c_{i}u_{k,m-\nu}
+ c_{0}u_{k,m-\nu-i} + \sum_{j=1}^{i-1} c_{j}u_{k,m-\nu-i+j} \right) \omega^{k}}
{u_{n-1,m-\nu}\omega^{n-1}+\cdots+u_{0,m-\nu}}.
$$

From \eqref{eq:ai}, this holds for $\nu=1$, allowing us to proceed by induction.

Thus
$$
\alpha_{i}^{(\nu)} -c_{i} = \frac{\sum_{k=0}^{n-1} \left( c_{0}u_{k,m-\nu-i} + \sum_{j=1}^{i-1} c_{j}u_{k,m-\nu-i+j} \right) \omega^{k}}
{u_{n-1,m-\nu}\omega^{n-1}+\cdots+u_{0,m-\nu}}.
$$

By Lemma~\ref{lem:2b}(b), for each $i=0,\ldots, n-1$, there exists at least
one $k \in \left\{ 0,\ldots, n-1 \right\}$ such that $u_{k,m-\nu-i} \geq 1$.
Thus both the denominator and numerator of this expression are positive for all
such $\nu$ and so $\alpha_{i}^{(\nu)}-c_{i}>0$.

From Lemma~\ref{lem:3},
\begin{equation}
\label{eq:rt-LB}
\omega>u_{n-1,m}t+c_{n-1} \geq (t+1)c_{n-1} \geq 2c_{n-1},
\end{equation}
since $m \geq n$, $t \geq 1$ and since $u_{n-1,m}\geq u_{n-1,n}=c_{n-1}$ by
Lemma~\ref{lem:2b}(a).

By Lemma~\ref{lem:2b}(d), we see that the coefficients of each $\omega^{k}$ in
the numerator of the expression above for $\alpha_{i}^{(\nu)}-c_{i}$ is at most
one more than the coefficient of the corresponding $\omega^{k}$ term in the
denominator. Furthermore, Lemma~\ref{lem:2b}(d) also allows us to show that
even in the worst case (when $i_{1}=n-1$ and we have equality in \eqref{eq:2b-d}
for all $0 \leq i \leq n-2$)
$$
\alpha_{i}^{(\nu)}-c_{i}
\leq 1+ \frac{-\omega^{n-1}+ \sum_{k=0}^{n-2} \omega^{k}}
{u_{n-1,m-\nu}\omega^{n-1}+\cdots+u_{0,m-\nu}}.
$$

So by \eqref{eq:rt-LB}, $\alpha_{i}^{(\nu)}-c_{i}<1$.

Using these expressions, we can write
$$
\frac{\alpha_{i}^{(\nu)}-a_{i}^{(\nu)}}{\alpha_{1}^{(\nu)}-a_{1}^{(\nu)}}
=\frac{\sum_{k=0}^{n-1} \left( u_{k,m-\nu-i} +\sum_{j=1}^{i-1} c_{j}u_{k,m-\nu-i+j} \right) \omega^{k}}
{\sum_{k=0}^{n-1} u_{k,m-\nu-1} \omega^{k}}
=\alpha_{i-1}^{(\nu+1)},
$$
for $i=2,\ldots,n-1$.

Lastly,
$$
\alpha_{n-1}^{(\nu+1)}=\frac{1}{\alpha_{1}^{(\nu)}-a_{1}^{(\nu)}}
=
\frac{u_{n-1,m-\nu}\omega^{n-1}+\cdots+u_{0,m-\nu}}
{\sum_{k=0}^{n-1} u_{k,m-(\nu+1)} \omega^{k}},
$$
from our above expression for $\alpha_{i}^{(\nu)}-c_{i}$. The desired
expression for $\alpha_{n-1}^{(\nu+1)}$ holds from the recurrence relation
\eqref{eq:recur} defining the $u_{k,m-\nu}$'s applied to the numerator.

\subsection{Induction from $\balpha^{(\nu)}$ to $\balpha^{(\nu+1)}$:
$m-n+1 \leq \nu \leq m-1$}
$\,$

From our previous induction, we know that for $\nu=m-n+1$,
$$
\alpha_{i}^{(\nu)}
= \left\{
\begin{array}{ll}
\omega^{-i}+\displaystyle\sum_{j=1}^{i} c_{j}\omega^{j-i}
& \mbox{if $1 \leq i \leq m-\nu$}, \\
\omega^{n-i}
- \displaystyle\sum_{j=i+1}^{n-1} \left( u_{j}t+c_{j} \right) \omega^{j-i}
& \mbox{if $m-\nu+1 \leq i \leq n-1$}.
\end{array}
\right.
$$

Now suppose for some $m-n+1 \leq \nu \leq m-1$, this holds. We proceed to
determine $\ba^{(\nu)}$ and $\balpha^{(\nu+1)}$.

From \eqref{eq:rt-LB},
$$
0 < \omega^{-i}+\sum_{j=1}^{i-1} c_{j}\omega^{j-i}
< 1,
$$
and so $a_{i}^{(\nu)}=c_{i}$ for $1 \leq i \leq m-\nu$. Thus, from Lemma~\ref{lem:3},
we have
$$
\ba^{(\nu)} = \left( c_{1},\ldots,c_{m-\nu}, u_{m-\nu+1,m}t+c_{m-\nu+1},
\ldots,u_{n-1,m}t+c_{n-1} \right).
$$

As before, we first consider $\alpha_{1}^{(\nu)}-a_{1}^{(\nu)}=1/\omega$. We
then verify that the above expressions for $\alpha_{i}^{(\nu)}$ hold for
$\alpha_{i}^{(\nu+1)}$ as well.

Therefore, the expressions above also hold for $\alpha_{i}^{(m)}$, which we see
equals $\alpha_{i}^{(0)}$, establishing our claim that this Jacobi-Perron
expansion is periodic with preperiod length $0$ and period length $m$.

\subsection{The Hasse-Bernstein unit}

From the expressions for the $\alpha_{i}^{(0)}$'s in Remark~\ref{rem:1}, they
are algebraic integers and so we can find the Hasse-Bernstein unit from the JPA
expansion. To do so, note that
\begin{align*}
\epsilon &=\prod_{\nu=0}^{m-1} \alpha_{n-1}^{(\nu)}
=\omega^{n-1} \prod_{\nu=1}^{m-n+1} \frac{u_{n-1,m+1-\nu}\omega^{n-1}+\cdots+u_{0,m+1-\nu}}
{u_{n-1,m-\nu}\omega^{n-1}+\cdots+u_{0,m-\nu}} \\
&=\omega^{n-1} \frac{u_{n-1,m}\omega^{n-1}+\cdots+u_{0,m}}
{u_{n-1,n-1}\omega^{n-1}+\cdots+u_{0,n-1}}
=u_{n-1,m}\omega^{n-1}+\cdots+u_{0,m},
\end{align*}
as stated.

\section{Remarks on Hasse-Bernstein Units}

\subsection{Hasse-Bernstein units for $n=3$}

For $n=3$, $f(X)$ in Theorem~\ref{thm:1} has negative discriminant if and only
if $c_{2} \geq \lceil c_{1}^{2}/4 \rceil$. In this case, the unit group in the
ring of integers of $\bbQ(\omega)$ is of rank $1$. Our calculations in PARI
suggest that, like the units in \cite{LR}, the Hasse-Bernstein unit is a
fundamental unit in the ring $\bbZ[\omega]$.

When we consider the full ring of integers of $\bbQ(\omega)$, it appears that
the Hasse-Bernstein unit is often a fundamental unit of the unit group.
However, it appears that for each tuple $\left( c_{0}=1, c_{1}, c_{2} \right)$
and each $m \geq 3$ satisfying the conditions of Theorem~\ref{thm:1}, there are
infinitely many positive integers $t$ such that the index is $2$.

E.g., for $m=6$ with $\left( c_{0}, c_{1}, c_{2} \right) = \left( 1,0,1 \right)$,
define the sequence: $a_{0}=-8, a_{1}=33$ and $a_k=-4a_{k-1}-a_{k-2}+4$ for $k \geq 2$.
If $t=12a_k^2 - a_{k-1}^2 - 14a_{k} + 2a_{k-1}+2$ for $k \geq 1$, then the group
generated by the Hasse-Bernstein unit has index $2$ in the full unit group.

Also for $m=3$ with $\left( c_{0}, c_{1}, c_{2} \right) = \left( 1,0,1 \right)$
and $t=k^{3}-1$ with $k \geq 2$, the index is $3$.

%

\subsection{Hasse-Bernstein units from \cite{LR}}
\label{subsect:units-LR}

Something similar happens with the Hasse-Bernstein units in \cite{LR}.

For example, with their family for their $4m+1$ where $m=1$, we define the
sequence: $a_{0}=0, a_{1}=-3$ and $a_{k}=6a_{k-1}-a_{k-2}-2$ for
$k \geq 2$. If $c=4a_{k}^2-a_{k}a_{k-1}-6a_{k}+a_{k-1}+2$, then the group
generated by the Hasse-Bernstein unit has index $2$ in the full unit group.

Also with their family for their $4m+1$ where $c=3$ and their $m$ is odd, the group
generated by the Hasse-Bernstein unit has index $2$.

For their family for their $3m+1$, we give here just two examples. With $c=2$
and $c=779$, both with their $m=1$, the groups generated by the Hasse-Bernstein
unit have index $3$ and $5$ respectively.

These examples answer negatively the questions in \cite{LR} about their
Hasse-Bernstein units being fundamental units in $\bbQ(\omega)$.

\section{Remarks on $f(X)$ in Theorem~\ref{thm:1}}

\subsection{Irreducibility of $f(X)$}

We have not shown that the polynomials, $f(X)$, in Theorem~\ref{thm:1}
are irreducible. Extensive computation suggests that they are, as
well as the following more general conjecture.

\begin{conjecture}
\label{conj:1}
For $n \geq 2$, let
$$
g(X)=X^{n}-a_{n-1}X^{n-1}-\cdots-a_{0} \in \bbZ[X],
$$
where $a_{n-1} \geq a_{n-2} \geq \cdots \geq a_{1} \geq 0$ and $a_{n-1} \geq a_{0} \geq 1$.
Either\\
{\rm (a)} if $a_{n-1}=a_{n-2}$, then $g(X)$ may be divisible by $\Phi_{m}$,
the cyclotomic polynomial of order $m$ where $m|n$ and $1<m<n$,\\
{\rm (b)} if $n$ is even, then $g(X)$ may be divisible by $\Phi_{2}(X)$,\\
{\rm (c)} if $n \equiv 5 \bmod 6$, then $g(X)$ may be divisible by $\Phi_{6}(X)$, or \\
{\rm (d)} $g(X)$ is irreducible in $\bbZ[X]$.
\end{conjecture}

%

\subsection{Roots of $f(X)$}

In most other work on periodic JPA expansions, $\balpha^{(0)}$ is associated
with either an $n$-th root of a rational number or a Pisot number. Here,
in many cases, $\omega$ is a Pisot number, but that is not always the case.
This is another novel aspect to our examples.

For example, if $n=4$, $\left( c_{0}, c_{1}, c_{2}, c_{3} \right) = (1,0,1,1)$,
$m=4$ and $t=1$, then $f(X)=X^{4} - 2X^{3} - 2X^{2} - 2$ also has a root at
$x=-1.13418\ldots$.

So there exist polynomials satisfying the conditions in
Conjecture~\ref{conj:1} with a real root less than $-1$ when $n$ is even. Such
roots appear to be bounded below by the negative real root of $x^{n-1}+x^{n-2}+1$.
However, for $n$ odd, such polynomials appear to have no real roots less than
$-1$.

There can also be non-real roots of $f(X)$ outside the unit circle. E.g., if $n=4$, $\left( c_{0}, c_{1}, c_{2}, c_{3} \right)
= (1,0,0,1)$,
$m=5$ and $t=10$, then $f(X)=X^{4} - 11X^{3} - 10X^{2} - 11$ has a pair of
non-real roots with absolute value $1.1908\ldots$.

\begin{center}
{\sc Acknowledgements}
\end{center}

I thank Claude Levesque for inspiring this work and for his stimulating
conversations on this subject, as well as the referee for their thoughtful
report and careful reading of the manuscript. Both played valuable, and much
appreciated, roles in this work, its results and its presentation.

\vspace{3.0mm}

\bibliographystyle{amsplain}

\end{document}